\documentstyle[12pt]{article}
\catcode`\@=11
\font\twelveeuf=eufm10 at 12pt
\font\teneuf=eufm10
\font\seveneuf=eufm7

\newfam\euffam
\textfont\euffam=\twelveeuf  \scriptfont\euffam=\teneuf
  \scriptscriptfont\euffam=\seveneuf
\def\Frak{\ifmmode\let\next\Frak@\else
  \def\next{\errmessage{Use \string\Frak\space only in math mode}}\fi\next}
\def\Frak@#1{{\Frak@@{#1}}}
\def\Frak@@#1{\fam\euffam#1}
\catcode`\@=12
\def\Bbb{\bf}
\def\hook{ \, \mbox{\raisebox{-.19ex}{\rule[.1mm]{2.15mm}
{.1mm}}\raisebox{-.15ex}{\rule{0.1mm}{2.3mm}}}~}
\newcommand{\be}[1]{\begin{equation} \label{#1}}

\def\ee{\end{equation}}
\def\bea{\begin{eqnarray*}}
\def\eea{\end{eqnarray*}}

\newtheorem{thm}{Theorem}[section]
\newtheorem{propn}[thm]{Proposition}
\newtheorem{cor}[thm]{Corollary}

\newtheorem{defn}{Definition}

\newenvironment{proof}{\medskip \noindent {\bf Proof.}}{\hfill \rule{.5em}{1em}
\\}
\def\+{\oplus}

\def\*{^{\ast}}

\def\Sg{\Sigma}
\def\bp{{\Bbb P}}

\def\bc{{\Bbb C}}

\def\J{{\cal J}}

\def\S{{\cal S}}
\begin{document}
\sloppy
\title{A K\"ahler Structure on the Space\\ of String World-Sheets}
\author{ Claude LeBrun\thanks{Supported
in part by  NSF grant DMS-92-04093.}\\Department of Mathematics\\SUNY Stony
 Brook}

\date{}
\maketitle

\begin{abstract} Let $(M,g)$ be an oriented
 Lorentzian 4-manifold, and consider the
space $\S$ of oriented, unparameterized time-like 2-surfaces in $M$
(string world-sheets) with
fixed boundary conditions.
Then the infinite-dimensional manifold
$\S$ carries a natural complex structure
and a compatible (positive-definite) K\"ahler metric $h$ on $\S$ determined
by the Lorentz metric $g$.
 Similar results are
proved for other dimensions and signatures, thus
generalizing results of Brylinski \cite{bryl} regarding knots
in 3-manifolds.
Generalizing the framework of Lempert \cite{lemp}, we also investigate
the precise sense in which $\S$ is
 an infinite-dimensional complex manifold.
\end{abstract}
\vfill

\begin{quote}
\begin{small}
{ Running title:}
{\sc The Space of String World-Sheets}
\end{small}
\end{quote}

\bigskip
 \pagebreak

\section{Introduction}
Given a collection of circles in
a 4-dimensional oriented Lorentzian space-time, one may consider the
space $\cal S$ of unparameterized oriented time-like compact 2-surfaces with
the given circles as boundary. The main purpose of
the present  note is to endow $\cal S$ with  the structure
of an infinite-dimensional K\"ahler manifold--- i.e.  with
both a  complex
structure and a Riemannian metric for which this
complex structure is covariantly constant.
This was  motivated by
a construction  of Brylinski \cite{bryl}, whereby a K\"ahler
structure is given to
 the space of knots in a Riemannian 3-manifold.
In fact, our  discussion will be structured so as to apply  to
codimension 2 submanifolds of a
  space-time of arbitrary
dimension and  metrics of arbitrary signature, with the
proviso that we  only consider those submanifolds
for which the normal bundle is orientable and has
 (positive- or negative-)definite
induced metric; thus Brylinski's construction becomes subsumed as a special
case.

As the reader will  therefore see,
complex manifold theory thus comes naturally into play when one
studies codimension 2 submanifolds of a space-time. On the other hand,
  complex manifold theory makes a quite different kind of
appearance when one attempts to study the intrinsic
geometry of 2-dimensional manifolds. If some interesting modification
of string theory could be found which invoked
both of  these observations simultaneously, one
might hope to thereby explain the puzzling four-dimensionality of the
observed world.

 Many of the key technical ideas in the present note  are
straightforward generalizations of arguments due to
 L\'aszl\'o  Lempert
\cite{lemp}, whose  lucid   study of Brylinski's complex structure
is based on the theory of twistor CR manifolds \cite{leb}.
One of the most striking features of the complex
structures in question is that, while they are formally integrable and
may even admit legions of  local holomorphic functions, they
do {\em not} admit enough finite-dimensional  complex
submanifolds to be locally modeled on any  complex
topological vector space.
This beautifully illustrates the fact, emphasized by  Lempert,
that the Newlander-Nirenberg Theorem \cite{nn} fails in
infinite dimensions.

\section{The Space of World-Sheets}
Let $(M, g)$ be a smooth oriented
pseudo-Riemannian n-manifold. We use the term {\em world-sheet} to refer
to a smooth compact
oriented codimension-2
submanifold-with-boundary
 $\Sg^{n-2}\subset M^n$ for which  the inner product induced
by $g$ on the conormal bundle
$$\nu^{\ast}_{\Sg}:=\{ \phi \in T^{\ast}M|_{\Sg}~~~|~~\phi|_{T\Sg}\equiv 0
\} $$
of $\Sg$
is definite at each point. If $g$ is Riemannian, this just means an
oriented submanifold of codimension 2;
on the other hand, if $(M,g)$ is a Lorentzian 4-manifold,
a world-sheet is exactly an oriented time-like 2-surface.

\begin{defn} Let $(M, g)$ be a smooth oriented
pseudo-Riemannian n-manifold, and let $B^{n-3}\subset M^n$ be a smooth
codimension-3 submanifold which is  compact, without boundary. We will
then let  $\S_{M,B}$ denote the space of smooth oriented
 world-sheets $\Sg^{n-2}\subset M$ such that
$\partial \Sg =B$.
\end{defn}

Of course, this space is sometimes empty--- as happens,  for example,  if $B$
is a single space-like circle in Minkowski 4-space. This said,
$\S_{M,B}$ is automatically a  Fr\'echet manifold, and its tangent space at
$\Sg$ is
$$T_{\Sg}\S_{M,B}=\{ v\in \Gamma (\Sg, C^{\infty}(\nu_{\Sg}))~~|~~
v|_{\partial \Sg}\equiv 0\}~ .$$
Indeed, if we choose a tubular neighborhood of $\Sg$ which is
identified with the normal bundle of an open extension
 $\Sg_{\varepsilon}$ of $\Sg$ beyond its boundary, every section
of $\nu_{\Sg}\to \Sg$ which vanishes on $\partial \Sg$ is thereby
identified with an imbedded submanifold of $M$, and this
submanifold is still a world-sheet provided the $C^1$ norm of the
section is sufficiently small. This provides $\S_{M,B}$ with
charts which take values in Fr\'echet spaces, thus giving it the desired
manifold structure.

Since the normal bundle
$\nu_{\Sg} =(\nu_{\Sg}^{\ast})^{\ast}=TM/T\Sg =(T\Sg)^{\perp}$ of  our
world-sheet is of rank 2 and comes equipped with
 an orientation as well as a metric induced by
$g$, we may identify $\nu_{\Sg}$ with a complex line bundle by
taking $J: \nu_{\Sg}\to\nu_{\Sg}$, $J^2=-1$ to be rotation by $+90^{\circ}$.
This
then defines an endomorphism $\J$ of $T\S$ by
$$\J : T_{\Sg}\S_{M,B}\to T_{\Sg}\S_{M,B}: v \to J\circ v ~ .$$
Clearly $\J^2=-1$, so that $\J$ gives $\S$ the structure of an almost-complex
Fr\'echet manifold--- i.e. every tangent space of $\S$ can now be thought of as
a complex Fr\'echet space by defining $\J$ to be multiplication by
$\sqrt{-1}$. In the next sections, we shall investigate the
integrability properties of this almost-complex structure.

\section{Integrability of the Complex Structure}
Let $(M,g)$ denote, as before, an oriented pseudo-Riemannian manifold.
Let $Gr_2^{+}(M)$
denote the bundle of oriented  2-planes in $T^{\ast}M$ on which
the inner product induced by $g$ is definite. This
smooth ($3n-4$)-dimensional manifold then has a natural CR structure
\cite{leb,ros} of codimension $n-2$. Let us review how this
comes about.

Let $\hat{N}\subset [(\bc\otimes T^{\ast}M)- T^{\ast}M]$
denote the set of non-real null covectors of $g$, and let $N\subset
 \bp (\bc\otimes T^{\ast}M)$ be its image in the fiber-wise projectivization
of the complexified cotangent bundle. There is then a natural identification
of ${N}$ with $Gr_2^{+}(M)$. Namely, using pairs $u,v\in T_x^{\ast}M$ of
real covectors
satisfying
 $\langle u, v\rangle =0$  and $\langle u, u\rangle =\langle v, v\rangle$,
 we define a bijection between these two spaces by
$$Gr_2^{+}(M) \ni
\mbox{oriented span}(u,v)\leftrightarrow [u+iv]\in {N}\subset \bp
(\bc\otimes T_x^{\ast}M)$$
which
is  independent of the representatives  $u$ and $v$. But, letting
$\vartheta=\sum p_jdx^j$ denote the canonical complex-valued 1-form
on the total space of $\bc\otimes T^{\ast}M\to M$, and letting $\omega$ be
the restriction
 of $d\vartheta$ to $\hat{N}$, the distribution
$$ \hat{D}=\ker (\omega:\bc\otimes T\hat{N}\to \bc\otimes T^{\ast}\hat{N})$$
is involutive by virtue of the fact that $\omega$ is closed; since
$ \hat{D}$ also contains no non-zero real vectors as a consequence of
the fact that $\hat{N}\cap T^{\ast}M=\emptyset$, $ \hat{D}$ is a CR structure
on $\hat{N}$, the codimension of which can be checked to be
$n-2$. This CR structure  is invariant under the natural action
of $\bc^{\times}$ on $\hat{N}$ by scalar multiplication, and
thus descends to a CR structure $D$ on $N=Gr_2^{+}(M)$, again of codimension
$n-2$. Moreover, $\vartheta|_{\hat{N}}$ descends to $N$ as a CR
line-bundle-valued 1-form
$$\theta\in \Gamma (N, {\cal E}^{1,0} (L))~,~~ \bar{\partial}_b\theta =0~,
$$
where, letting $T^{1,0}N:=({\Bbb C}\otimes TN)/D$,
 $L^{\otimes (n-1)}=\bigwedge^{(2n-3)}T^{1,0}N$,
${\cal E}^{1,0} (L):= C^{\infty} (L\otimes (T^{1,0}N)^{\ast})$,
and $\bar{\partial}_b$ is naturally induced by $d|_D$.

The CR structure  $D$ of $N$  may be  expressed in the form
$$D=\{ v-iJv~| ~~v\in H\}$$
for a unique rank $2n-2$ sub-bundle
$H$ of the real tangent bundle $TN$ and a unique endomorphism $J$ of $H$
satisfying $J^2=-1$.
In these terms the geometric
 meaning of the  CR structure of $N$ is
fairly easy to describe. Indeed, if   $\varpi : Gr_2^{+}(M)\to M$
is the tautological projection, then
$H_P = (\varpi_{\ast P})^{-1}(P)$ for every  oriented definite 2-plane
$P\subset TM$.
On vertical vectors,  $J$ acts  by the
standard complex structure on the quadric fibers of $N\to M$;
whereas $J$ acts on horizontal vectors by $90^{\circ}$ rotation in the
2-plane $P\subset TM$. This point of view, however,
obscures the fact that $D$ is both involutive
and
unaltered by  conformal changes $g\mapsto e^fg$.

A compact  $(n-2)$-dimensional submanifold-with-boundary  $S\subset N$,
 will be called a  {\em transverse sheet} if its tangent space
is everywhere tansverse to the CR tangent space of $N$:
$$TN|_S= TY\oplus H|_Y~ .$$ As before,
let $B^{n-3}\subset M$ denote a  compact codimension-3
submanifold, and let $\varpi :N\to M$ be the canonical projection.
We will then let
$\hat{\S }_{N,B}$ denote the set of
transverse sheets $S\subset N$
such that $\varpi$ maps $\partial S$ diffeomorphically onto
$B$. Thus $\hat{\S }_{N,B}$ is a Fr\'echet manifold whose
tangent space at $S$ is given by
$$T\hat{\S }_{N,B}|_S=\{ v\in \Gamma (S, C^{\infty}(H|_S))~
 |~\varpi_{\ast}(v|_{\partial S})\equiv 0\} ~ ,$$
and hence $J: H\to H$ induces an almost-complex structure
$\hat{\cal J}$ on
$ \hat{\S }_{N,B}$ by $\hat{\cal J} (v):= J\circ v$.

\begin{propn} The almost-complex structure $\hat{\cal J}$ on the space
$ \hat{\S }_{N,B}$ of transverse sheets is formally integrable--- i.e.
$$\tau (v,w):= \hat{\cal J}[v,w]
-[v,\hat{\cal J}w] -[\hat{\cal J}v,w]
-\hat{\cal J} [\hat{\cal J}v,\hat{\cal J}w] =0 $$
for all smooth  vector fields $v$, $w$ on $\hat{\S }_{N,B}$.
\end{propn}
\begin{proof}
The Fr\"ohlicher-Nijenhuis torsion $\tau (v,w)$ is tensorial in the sense
that its value at $S$ only depends on the values of $v$ and $w$ at
$S$.  Given $v_S,w_S \in \{ v\in \Gamma (S, C^{\infty}(H|_S))~
 |~\varpi_{\ast}(v|_{\partial S})\equiv 0\}$, we will now
define
preferred extensions of them as vector fields near $S\in  \hat{\S }_{N,B}$
in such a manner as to simplify the computation of
$\tau (v,w)= \tau (v_S,w_S)$.
     To do this, we may first use a partition of unity to
extend $v_S$ and $w_S$ as sections
$\hat{v}, \hat{w}\in \Gamma (N, C^{\infty}(H))$ defined on all of
of $N$ in such a manner  that  $\hat{v}$ and $\hat{w}$
are tangent to the fibers of $\varpi$ along all of $\varpi^{-1}(B)$.
Now  let $U\subset N$ be a tubular neighborhood
of $S$ which is  identified with the normal bundle $H$ of
some open extension $S_{\epsilon }$ of $S$, and
let $\hat{U}\subset \hat{\S }_{N,B}$ be the
set of transverse sheets $S'\subset U$.
We may now define our preferred extensions  of
$v$ and $w$ of $v_S$ and $w_S$ on the domain $\hat{U}$ by letting
the values of $v$ and $w$ at $S'\subset U$
be the restrictions of $\hat{v}$ and $\hat{w}$
 to ${S'}$. Notice that $[ v,w ]$ is then precisely the vector field
on $\hat{U}$ induced by $[\hat{v}, \hat{w}]$, whereas
$\hat{\cal J}v$ is the vector field induced by $J\hat{v}$.
Since the integrability condition for $(N,D)$ says that
$$J([\hat{v},\hat{w}]  -[J\hat{v},J\hat{w}])
=[\hat{v},J\hat{w}] +[J\hat{v},\hat{w}]~ , $$
it therefore follows that
$$\hat{\cal J}[v,w]
-\hat{\cal J} [\hat{\cal J}v,\hat{\cal J}w] =
[v,\hat{\cal J}w] +[\hat{\cal J}v,w] ~ ,$$
so that $\tau (v,w)=0$, as claimed.
\end{proof}

We now observe that there is a canonical imbedding
\bea {\S }_{M,B}&\stackrel{\Psi}{\hookrightarrow} &\hat{\S }_{N,B}\\
\Sg &\mapsto & \nu_{\Sg }
\eea
obtained by sending a world-sheet to its normal-bundle,
thought of as the image of a  section of $Gr_2^+(M)|_{\Sg }= N|_{\Sg }$;
thought of in this way, it is easy to see that $\nu_{\Sg}\subset N$
is a transverse submanifold.

\begin{thm} The imbedding $\Psi$ realizes $({\S }_{M,B}, {\cal J})$
as a complex submanifold of $(\hat{\S }_{N,B}, \hat{\cal J})$.
In particular, the almost-complex structure  ${\cal J}$
of ${\S }_{M,B}$ is formally integrable. \label{imb}
\end{thm}
\begin{proof}
The projection $\varpi: N\to M$ induces a map  $\hat{\varpi}:
\hat{\S }_{N,B}\to{\S }_{M,B}$
which is  a left inverse of $\Psi$ and satisfies $\hat{\varpi}_{ *}\hat{\cal J}
={\cal J}\hat{\varpi}_{ *}$. It therefore suffices to show that
the tangent space of the image of
$\Psi$  is  $\hat{\cal J}$-invariant.
Now the condition for a transverse sheet $S\subset N$
to be the $\Psi$-lifting of the
world-sheet  $\varpi(S)\subset M$ is exactly that  $\theta|_S\equiv 0$.
When $S$ satisfies this condition,  a connecting field
$v\in\Gamma (S, {\cal E}(H))$
then represents a vector $\hat{v}\in T\hat{\cal S}$  which is
 tangent to the image of $\Psi$ iff
\be (v\hook d\theta)|_{TS}+d(v\hook \theta)|_{TS}\equiv 0~ ;\label{leg} \ee
 the exterior derivative of $\theta$ may here be calculated in any
local trivialization for the line bundle $L$, since the
left-hand side rescales properly under changes of trivialization so as
define an $L$-valued 1-form on $S$. But since $\theta\in \Gamma (N, {\cal
E}^{1,0}
(L))$ satisfies $\bar{\partial}_b \theta =0$, it follows that
$$(Jv\hook d\theta)|_{TS}+d(Jv\hook \theta)|_{TS}=i(Jv\hook d\theta)|_{TS}+
id(Jv\hook \theta)|_{TS} $$
because $\theta$ and $d\theta$ are of types (1,0) and (2,0), respectively.
The tangent space of the image of $\Psi$
is therefore $\hat{\cal J}$-invariant, and the claim follows.
\end{proof}

\begin{defn} Let $({\Frak X}, {\Frak J})$
be an almost-complex Fr\'echet manifold, and let
$f: U\to {\Bbb C}$ be a differentiable function defined on an
open subset of ${\Frak X}$. We will say that
$f$ is ${\Frak J}$-holomorphic if
$$({\Frak J}v)f=ivf~~\forall v\in TU~.$$
\end{defn}

\begin{defn}
An almost-complex Fr\'echet manifold
$({\Frak X}, {\Frak J})$ is
called {\em weakly integrable} if for each real tangent vector
$w\in T{\Frak X}$ there is a ${\Frak J}$-holomorphic
 function $f$ defined on a neighborhood
of the base-point of $w$ such that $wf\neq 0$.
\end{defn}

\begin{thm} Suppose that $(M,g)$ is real-analytic. Then
$(\hat{\S }_{N,B}, \hat{\cal J})$ is weakly integrable.\label{wint}
\end{thm}
\begin{proof}
If  $(M,g)$ is real-analytic, so is the CR manifold $(N,D)$,
and we can therefore
realize $(N,D)$
as a real submanifold of a complex manifold $(2n-3)$-manifold
$\cal N$. This can even be done explicitly by
taking $\cal N$ to be a space of  complex null geodesics
for a suitable complexification of $(M,g)$.

Now let $S\subset N\subset {\cal N}$ be any transverse sheet.
Then there is a neighborhood $V\subset {\cal N}$ of $S$ which can
be holomorphically imbedded in some ${\Bbb C}^{\ell}$.
Indeed, let $Y\subset {\cal N}$ be  a totally real $(2n-3)$-manifold
containing $S$, let $f: Y\to  {\Bbb R}^{\ell}$ be
a smooth imbedding, and let $Y_0$ be a precompact
neighborhood of $S\subset Y$ with smooth boundary.
 By \cite{wel},the component functions  $f^j|_{Y_0}$ are  limits
in the $C^1$ topology
of the restrictions of holomorphic functions. Using such an approximation
of $f$,
 we may therefore  imbed $Y_0$ as a totally real submanifold
of  ${\Bbb C}^{\ell}$ by a map which
extends   holomorphically to a neighborhood of
$Y_0$, and this holomorphic extension then automatically
yields a holomorphic imbedding
of some open neighborhood $V\supset S$ in ${\Bbb C}^{\ell}$.

Now suppose that  $v$ is a smooth section of $H$ along $S$.
We may express $v$ uniquely  as $u+Jw$, where $u$ and $w$ are
tangent to the $Y$. By changing $Y$ if necessary, we can furthermore
assume that $u\not\equiv 0$.
Let $F$ be a smooth function on $Y$ which
vanishes on $S$ and such that the derivative $vF$ is
non-negative and supported near some interior point of $S$;
 and let $\varphi$ be a real-valued smooth $(n-2)$-form
on $Y$ whose restriction to $S$ is positive on the
support of $vF$. Set $\psi =F\varphi$. Using
our imbedding of $Y$ in ${\Bbb C}^{\ell}$,
we can express $\psi$ as a family of component functions---
e.g. by arbitrarily declaring that all contractions
of $\psi$ with  elements of the normal bundle
$(TY)^{\perp}\subset T{\Bbb C}^{\ell }$ shall vanish.
But, again by \cite{wel}, these component functions are $C^1$-limits
on $Y_0$ of restrictions of
holomorphic functions from a neighborhood of $Y_0\subset {\Bbb C}^{\ell}$.
Thus, by perhaps replacing $V$ with a smaller neighborhood,
there is a holomorphic $(n-2)$-form $\beta$ on $V$ which approximates
  $\psi$ well enough that $$\Re e \int_Sv\hook d\beta > \frac{1}{2}
\int_Su\hook d\psi >0$$ and $$\Re e
\int_{\partial S}v\hook \beta > - \frac{1}{2}
\int_Su\hook d\psi~ .$$

Let $\hat{V}:=\{S'\in {\cal S}_{N,B}~|~S'\subset V\}$,
and define  $f_{\beta}: \hat{V}\to {\Bbb C}$
by $f_{\beta}(S')=\int_{S'}\beta$. Then $f_{\beta}$
is a holomorphic function on the open set
$\hat{V}\subset {\cal S}_{N,B}$. Indeed, if
$\gamma$ is {\em any} smooth $(n-2)$-form on $V$, and if we set
$f_{\gamma}(S')=\int_{S'}\gamma$, then, for $S'\subset V$,
the derivative of $f_{\gamma}$ in the direction of
$w\in \Gamma (S', C^{\infty}(H))$, $\varpi_*(w)|_{\partial S'}\equiv 0$,
 is given by
$$wf_{\gamma}|_S= \int_{S}w \hook d\gamma +
\int_{\partial S'}w \hook\gamma ~ ;$$
for if $w$ is extended to $V$ as a smooth
vector field $\hat{w}$ tangent  to the
fibers of $\varpi$ and $S_t$ is obtained by pushing $S'$
along the flow of the vector field $\hat{w}$, then
\begin{eqnarray*} wf_{\gamma}|_{S'}&=&
 \left.\frac{d}{dt}\left[\int_{S_t}\gamma\right]\right|_{t=0}
\\&=& \int_{S'}\pounds_{\hat{w}} \gamma
\\&=& \int_{S'}[\hat{w} \hook d\gamma +d(\hat{w} \hook\gamma )]
\\&=& \int_{S'}w \hook d\gamma +\int_{\partial S'}w \hook\gamma ~ .
\end{eqnarray*}
But since $\beta$ is the restriction of a
 holomorphic $(n-2)$-form from a region of $\cal N$,
it therefore follows that
\begin{eqnarray*}(\hat{\cal J}w)f|_{S'}&=&
 \int_{S'}(Jw )\hook d\beta +\int_{\partial S'}(Jw) \hook\beta
\\&=& i\int_{S'}w \hook d\beta +i\int_{\partial S'}w \hook\beta
\\&=& iwf|_{S'}~ ,\end{eqnarray*}
showing that the function $f_{\beta}$ induced by $\beta$  is
 $\hat{\cal J}$-holomorphic, as claimed.

However, we have also carefully chosen $\beta$ so that the real part
of the expression $\int_{S}v \hook d\beta +\int_{\partial S}v \hook\beta
= vf_{\beta}$ is positive.
For every real tangent vector $v$ on ${\cal S}_{N,B}$,
one can thus find a locally-defined
 $\cal J$-holomorphic function whose derivative
is non-trivial in the  direction $v$. Hence  $\hat{\cal J}$
is weakly integrable, as claimed.
 \end{proof}

\begin{cor} If $(M,g)$ is real-analytic, then
$({\S }_{M,B}, {\cal J})$ is weakly integrable.
\end{cor}
\begin{proof} By
   Theorem \ref{imb} and \ref{wint}, $({\S }_{M,B}, {\cal J})$
can be imbedded
in the weakly integrable almost-complex manifold
$(\hat{\S }_{N,B}, \hat{\cal J})$.
Since the restriction of a holomorphic function to an almost-complex
submanifold is holomorphic, it follows that $({\S }_{M,B}, {\cal J})$
is weakly integrable.
 \end{proof}

One might instead ask whether $({\S }_{M,B}, {\cal J})$ is
{\em strongly integrable}--- i.e. locally biholomorphic to a ball in some
complex vector space. The answer is {\bf no};   in contrast to
any strongly integrable almost-complex manifold,
 $({\S }_{M,B}, {\cal J})$ contains  very few
finite-dimensional complex submanifolds:

\begin{propn} Suppose that $(M,g)$ is real-analytic.
At a generic point $S\in {\S }_{M,B}$, a generic $(n-1)$-plane
is not tangent to any $(n-1)$-dimensional ${\cal J}$-complex submanifold.
\end{propn}
\begin{proof}
Let $S\subset M$ be a world-sheet which is {\sl not} real-analytic near
$p\in S$. Let $q\in N=Gr_2^+(M)$ be given by $q=T_pS^{\perp}$, and let
$v_1, \ldots , v_{n-1}\in H_q$ be a set of real vectors such that the
$v_j+iJv_j$ form a basis for $D_q$. Extend these vectors as sections
$\hat{v}_j$ of $H|_{\Psi(S)}$ which  satisfy equation (\ref{leg})
along the sheet; this may be done, for example, by first extending each
$v_j$ to just a 1-jet at $p$ satisfying  (\ref{leg}) at $p$, projecting this
to a 1-jet via $\varpi$ to yield
a 1-jet of a normal vector field on $S\subset M$, extending this 1-jet as a
section of the normal bundle of $S$, and finally lifting this section
using $\Psi_{*}$. Let $u_1, \ldots , u_n$ be the
elements of  $T_S{\S }_{M,B}$ represented by
$\hat{v}_1, \ldots , \hat{v}_{n-1}$,
and let $P\subset T^{1,0}_S{\S }_{M,B}$
be  spanned by $u_1-i{\cal J}u_1, \ldots , u_n-i{\cal J} u_n$.

Now suppose there were an $\cal J$-holomorphic submanifold
$X\subset {\S }_{M,B}$
 through $S$ with
(1,0)-tangent plane equal to $P$. Then $X$  represents a
family of transverse sheets in $N$ which foliates  a neighborhood $U\subset
N$ of $q$;
moreover, because $X$  represents a {\sl holomorphic}  family,
the  leaf-space projection $\ell : U\to X$ is CR in the sense that
$\ell_*(D)\subset T^{0,1}X$. Since we have assumed that
$(M,g)$ is real-analytic, we may also assume that $U$ has
a  real-analytic CR imbedding $U\hookrightarrow {\Bbb  C}^{2n-3}$.
Moreover the twistor CR manifold $N$ is automatically
``anticlastic,'' by which I mean that the Levi form
${\cal L}: D\to TN/H: v\mapsto  [v,\bar{v}] \bmod H$
is {\sl surjective} at each point of $N$.
This gives rise to a Bochner-Hartogs
extension phenomenon:
every CR function on $U$ extends to a holomorphic function on some
neighborhood of $U\subset {\Bbb  C}^{2n-3}$. In particular, every
CR map defined on $U$ must be real-analytic, and this applies in particular
to the leaf-space projection $\ell$. Thus $\Psi(S)$ is real-analytic near
$q$, and $S$ is therefore real-analytic near $p$.
This proves the result by contradiction.
\end{proof}

\section{K\"ahler Structure}
The complex structure $\cal J$ on $\S$ depends only
on the conformal class $[ g]= \{ e^fg\}$ of our
metric, but we will now specialize by fixing a specific
 pseudo-Riemannian metric $g$. Our reason for
 doing so is that we thereby induce an
$L^2$-metric  on $\S$. Indeed,
each tangent space
$$T_{\Sg}\S_{M,B}=\{ v\in \Gamma (\Sg, C^{\infty}(\nu_{\Sg}))~~|~~
v|_{\partial \Sg}\equiv 0\}~$$ may be equipped with
 a positive-definite
 inner product by setting
$$h( v, w ) := \int_{\Sg}g(v, w) ~d\mbox{vol}_{g|_{\Sg }} .$$
We shall now see that this metric has some quite remarkable
properties.

\begin{thm} The Riemannian metric $h$ on $\S$ is  {\em K\"ahler}
 with respect to the
previously-defined complex structure ${\cal J}$.\end{thm}
\begin{proof}
Let
$\Omega$ denote the volume n-form of $g$, and
define a
2-form $\omega$ on $\S$  by
$$\omega  ( v, w ) := \int_{\Sg}(v\wedge w)\hook \Omega ~ .$$
Obviously,  $\omega$ is ${\cal J}$-invariant and
$$h( v, w )=\omega ({\cal J}v, w).$$
We therefore just need to show that $\omega$ is closed.\footnote{
The reader may ask whether it is actually legitimate to
call a Riemannian
manifold K\"ahler when the almost-complex structure in question is
at best weakly integrable. However, formal integrability and the
closure of the K\"ahler form are the only conditions
 necessary  to insure that
the  almost-complex structure tensor is parallel,
even in infinite dimensions.}
     To check this,
let us introduce the universal family
\setlength{\unitlength}{1ex}
\begin{center}\begin{picture}(20,17)(0,3)
\put(9.5,17){\makebox(0,0){${\cal F}$}}
\put(1.5,5){\makebox(0,0){${\cal S}$}}
\put(17,5){\makebox(0,0){$M$}}
\put(14.5,12){\makebox(0,0){$p$}}
\put(4.5,12){\makebox(0,0){${\pi }$}}
\put(10.5,15){\vector(2,-3){5.5}}
\put(8,15){\vector(-2,-3){5.5}}
\end{picture}\end{center}
 where the
fiber of $\pi$ over $\Sg \in \S$ is defined to be $\Sg \subset M$.
 We can then pull $\Omega$
back to ${\cal  F}$
to obtain a closed n-form $\alpha=p^{\ast}\Omega$ which
 vanishes on the boundary $B=\partial \Sg$ of every  fiber
of $\pi $. But $\omega$ is just obtained from $\alpha$ by integrating
on the fibers of $\pi $:
$$\omega=\pi _{\ast}\alpha~ .$$
Since $\pi _{\ast}$ commutes with $d$ on forms which vanish along the
fiber-wise boundary (cf. \cite{bt}, Prop. 6.14.1), it follows that
$$d\omega=d(\pi _{\ast}\alpha )= \pi _{\ast}d\alpha=
\pi _{\ast}d (p^{\ast}\Omega )=\pi _{\ast}p^{\ast}d\Omega=
\pi _{\ast}p^{\ast}0=0.$$
Thus $h$ is a K\"ahler metric, with K\"ahler form $\omega$.
\end{proof}

    To conclude this note, we now observe that $({\cal S}, h)$ is formally
of Hodge type---  but  non-compact, of course!

\begin{propn} Modulo  a  multiplicative constant, the
K\"ahler form $\omega$  of $h$
 represents an integer  class in cohomology.  If,
moreover,  $M$  is  non-compact, $\omega$ is actualy an exact form, and
its cohomology class thus vanishes.
\end{propn}
\begin{proof}
If $M$ is compact, we may assume that $g$ has total volume 1, so that
its volume form $\Omega$ then  represents an element of integer cohomology.
Since $\omega=\pi _{\ast}p^{\ast}\Omega$, its
cohomology class $[\omega ]=\pi _{\ast}p^{\ast}[ \Omega ]$
is therefore integral. If, on the other hand, $M$ is
non-compact, $\Omega =d\Upsilon$ for some
$(n-1)$-form $\Upsilon$, and hence $\omega = d(\pi _{\ast}p^{\ast}\Upsilon)$.
\end{proof}

\bigskip
\noindent
{\bf Acknowledgements.} The author would like to
thank L\'aszl\'o Lempert and
Edward Witten for their suggestions and encouragement.

\end{document}